\documentclass[4pt,runningheads,a4paper]{article}
\usepackage{graphicx}
\usepackage{amsmath}
\usepackage{amssymb}

\makeatletter

\newcommand{\Rmnum}[1]{\expandafter\@slowromancap\romannumeral #1@}
\makeatother

\usepackage{amsmath,amssymb,graphics,float}
\usepackage{pstricks,pst-node,pst-tree,pst-plot}
\usepackage{graphicx,epsfig,subfigure}

\usepackage{color}
\usepackage{graphicx}
\usepackage{amsmath}
\usepackage{amsfonts,amssymb}
\usepackage{upgreek}

\renewcommand{\tau}{\uptau}
\renewcommand{\omega}{\upomega}

\newcommand{\Id}{M}

\newcommand{\ww}{\hat{w}}

\newcommand{\tn}{\tau_{n}}
\newcommand{\spV}[1]{\left< #1 \right>}
\newcommand{\spH}[1]{\left( #1 \right)_{H}}

\newcommand{\Ih}{\hat{I}}

\renewcommand{\th}{\hat{t}}
\newcommand{\phih}{\hat{\phi}}

\newcommand{\psih}{\hat{\psi}}

\newcommand{\R}{{\mathbb R}}   
\renewcommand{\P}{{\mathbb P}}   

\usepackage{authblk}
\author[1]{M. A. Qureshi}
\author[2]{S. Hussain}
\author[3]{Ghulam Shabbir}
\affil[1,3]{Faculty of Engineering Sciences, GIK Institute of Engineering Sciences and Technology, Topi, Khyber Pakhtunkhwa, Pakistan. }
\affil[2]{Department of Mathematics, Mohammad Ali Jinnah University, Islamabad, Pakistan.
Emails: amergikian@yahoo.com 
}
\date{}
\title{Conservation of Hamiltonian using Continuous Galerkin Petrov time discretization scheme}
\begin{document}

\maketitle

\begin{abstract}
\noindent Continuous Galerkin Petrov time discretization scheme is tested on some Hamiltonian systems including simple harmonic oscillator, Kepler's problem with
different eccentricities and molecular dynamics problem. In particular, we implement the fourth order Continuous Galerkin Petrov time
discretization scheme and analyze numerically, the efficiency and conservation of Hamiltonian. A numerical comparison with some symplectic methods including Gauss
implicit Runge-Kutta method and general linear method of same order is given for
these systems. It is shown that the above mentioned scheme, not only preserves Hamiltonian but also uses the least CPU time compared with upto-date and optimized methods.
\end{abstract}

\centerline{{\bf Mathematics Subject Classification:} }
\vspace{.3cm} {\bf Keywords:} Hamiltonian systems, Continuous Galerkin Petrov time discretization, G-symplectic general linear methods, Runge-Kutta Mathod,  Simple harmonic oscillator, Kepler's problem and Molecular dynamics problem

\section{Introduction}
\label{intro}
Non-dissipative phenomena arising in the fields of classical mechanics, molecular dynamics, accelerator physics, chemistry and other sciences
are modeled by Hamiltonian systems.
Hamiltonian systems define equations of motion based on generalised co-ordinates $q_{i}=(q_{1},q_{2},\cdots,q_{n})$ and generalised momenta  $p_{i}=(p_{1},p_{2},\cdots,p_{n})$ and are given as,
 \begin{equation} \label{hamiltonian}
 \frac{dp_{i}}{dt}=-\frac{\partial{H}}{\partial{q_{i}}},\hspace{0.5in}\frac{dq_{i}}{dt}=\frac{\partial{H}}{\partial{p_{i}}},\hspace{1in}  i=1,\cdots,n,
\end{equation}
having $n$ degrees of freedom. $H:\mathbb{R}^{2n}\times \mathbb{R}^{2n}\to \mathbb{R}$ is the total energy of the Hamiltonian system. A separable Hamiltonian has
the structure
\begin{equation*}
 H(p,q)=T(p)+V(q)
\end{equation*}

in mechanics, $T=\frac{1}{2}P^T M^{-1} P$ represents the kinetic energy and $V$ being the potential energy. The Hamiltonian system in partitioned form takes the form
\begin{equation*} \label{hamiltonian_partitioned}
 \frac{dp_{i}}{dt}=-\triangledown_q V,\hspace{0.5in}\frac{dq_{i}}{dt}=\triangledown_p T = M^{-1} p.
\end{equation*}

The first observation is that, for autonomous Hamiltonian systems, $H$ is an invariant, thus by differentiating $H(p,q)$ with respect to time we have,
\begin{equation*}\label{energy_conservation}
 \frac{dH}{dt} = \displaystyle\sum_{i=1}^{n}\Big{(}\frac{\partial{H}}{\partial{p_{i}}}\frac{dp_{i}}{dt} + \frac{\partial{H}}{\partial{q_{i}}}\frac{dq_{i}}{dt}\Big{)} = 0 .
\end{equation*}
We can write $ y = (p,q)$, then \eqref{hamiltonian} can be written as,
\begin{equation*} \label{hamiltonian-compact}
y' = J^{-1} \nabla H,
\end{equation*}
where $'$ represents the derivative with respect to time, $\nabla$ is a gradient operator and $J$ is a skew symmetric matrix consisting of zero matrix $0$ and $n \times n$  identity matrix $I$,
\begin{equation*}\label{J_matrix}
J = \left[ \begin{array}{cc}
0 & I \\
-I & 0 \end{array} \right].
\end{equation*}
Another property of Hamiltonian systems is that its flow is symplectic, i.e. for a linear transformation $ \Psi:\mathbb{R}^{2n} \mapsto
\mathbb{R}^{2n}$, the jacobian matrix $\Psi'(y)$ satisfies
\begin{align*}
\Psi'^{T}(y) \mathbf{\mathrm{J}} \Psi'(y) =\mathbf{\mathrm{J}}.
\end{align*}

Conservation laws for Hamiltonian systems are generally lost while integrating these system. It is generally desirable to preserve the
underlying qualitative property of solutions of Hamiltonian systems. This is achieved by using symplectic integrators from the class of
one step, multistep and general linear methods. A lot of attention has been paid on the construction and implementation of such integrators,
for details see \cite{hairer_gni}, \cite{rf1}, \cite{rf2} and \cite{sanzserna}.

The continuous Galerkin Petrov time discretization scheme (cGP) was investigated in~\cite{Schieweck2010} for the system of ordinary differential equations (ODEs). In~\cite{Hussain2011}, this scheme was studied for the heat equation. In particular, the cGP(2) scheme has found to be 4th order accurate in the discrete time point and is A-stable method.

The objective of this paper is to provide analysis of cGP(2) scheme \cite{Schieweck2010,Hussain2011,Aziz1989,Vidar2006} on some Hamiltonian systems and comparing it with other symplectic methods of order four
including Gauss implicit Runge-Kutta method represented as irk4 \cite{butcher-book} and a g-symplectic general linear method represented by glm4 of same order developed in  \cite{butcher_yousaf} and \cite{yousafthesis}.
In section two a brief introduction about the methods is given. The tested problems of Hamiltonian systems along with numerical experiments of these methods
 on Hamiltonian systems are described in third section. Conclusion based on
numerical comparison of third section is given in fourth section.

\section{The Methods}

\subsection{{\em Continuous} Galerkin-Petrov method (cGP)
\label{cGP_method}
}
As a model problem we consider the ODE system given in~\eqref{hamiltonian}: {\em Find $u:[0,T_m]\to W$ such that}
\begin{equation} \label{cont_eq}
  \begin{array}{rcll}
    d_t u(t) &=& F(t, u(t)) \quad \text{for} \quad  t\in(0,T_m) , \\
    u(0) &=& 0
  \end{array}
\end{equation}
The {\em weak formulation} of problem \eqref{cont_eq} reads:
Find $u\in X$ such that $u(0)=u_0$ and
\begin{equation} \label{weak_eq}
 \int_0^{T_m} \spV{d_t u(t), v(t)} dt = \int_0^{T_m} \spV{F(t, u(t)), v(t)} dt
 \qquad\forall\;  v\in Y ,
\end{equation}
where $X$ denotes the  {\em solution space} and Y the test space.
To describe the time discretization of problem~\eqref{cont_eq} let us introduce the following notation.
We denote by $I = [0, T_m]$ the time interval with some positive final time $T_m$.
We start by decomposing the time interval $I$ into $N$ subintervals
$I_n := (t_{n-1},t_n)$, where $n\in\{1, \dots, N\}$ and
$$
 0=t_0 < t_1 < \dots < t_{N-1} <  t_N = T_m.
$$
In our time discretization, we approximate the {\em continuous solution}
$u(t)$ of problem~\eqref{cont_eq} on each time interval $I_n$ by a polynomial
function:
\begin{equation} \label{poly_ansatz}
   u(t) \approx u_{h}(t) := \sum_{j=0}^{k} U_n^j \phi_{n,j}(t)
   \qquad\forall \; t\in I_n ,
\end{equation}
where the ''coefficients'' $U_n^j$ are elements of the Hilbert space $W$
and the basis functions $\phi_{n,j} \in\P_{k}(I_n)$ are linearly independent
elements of the standard space of polynomials on the interval $I_n$ with
a degree not larger than a given order $k$.

For a given time interval $J\subset\R$ and a Banach space $B$, we introduce
the linear space of $B$-valued time polynomials with degree of at most $k$  as
$$
  \P_k(J,B) :=
  \left\{
    u: J\to B \,:\;
    u(t) = \sum_{j=0}^k U^j t^j\,,\;\forall\, t\in J,\; U^j\in B,
    \;\forall\, j
  \right\}.
$$
Now, the {\em discrete solution space} for the global approximation $u_h:I\to W$
is the space $X_h^{k}\subset X$ defined as
\begin{equation*} \label{Xt_def}
  X_h^{k} := \{ u\in C(I,W) :\; u\big|_{\bar{I}_n} \in \P_k(\bar{I}_n, W)
                 \quad\forall\; n=1,\dots,N \}
\end{equation*}
and the {\em discrete test space} is the space $Y_h^{k}\subset Y$ given by
\begin{equation*} \label{Yt_def}
  Y_h^{k} := \{ u\in L^2(I,W) :\; u\big|_{I_n} \in \P_{k-1}(I_n, W)
                 \quad\forall\; n=1,\dots,N \}.
\end{equation*}
The symbol $h$ denotes the {\em discretization parameter} which acts in the
error estimates as the maximum time step size
$\; h := \max_{1\le n\le N} h_n$, where $\;h_n := t_n - t_{n-1}$ is the
length of the $n$-th time interval $I_n$.

Let us denote by $X_{h,0}^{k} := X_h^{k}\cap X_0$ the subspace of $X_h^{k}$ with zero
initial condition.
Then, it is easy to see that the dimensions of the spaces $X_{h,0}^{k}$ and
$Y_h^{k}$ coincide such that it makes sense to consider the following
{\em discontinuous Galerkin-Petrov discretization} of {\em order $k$}
for the weak problem \eqref{weak_eq} :\;
Find $u_h\in u_0 + X_{h,0}^{k}$ such that
\begin{equation} \label{cGP_form}
 \int_0^{T_m} \spV{d_t u_h(t), v_h(t)} dt = \int_0^{T_m} \spV{F(t, u_h(t)), v_h(t)} dt
 \qquad\forall\;  v_h\in Y_h^{k}.
\end{equation}
We will denote this discretization as the ''{\em exact} cGP(k)-{\em method}''.
Since the discrete test space $Y_h^{k}$ is discontinuous, problem \eqref{cGP_form}
can be solved in a time marching process.
Therefore, we choose test functions
$v_h(t) = v \psi_{n,i}(t)$  with an arbitrary $v\in W$  and a scalar function
$\psi_{n,i} : I\to\R$ which is zero on $I\setminus \bar{I}_n$ and
a polynomial $\psi_{n,i}\in\P_{k-1}(\bar{I}_n)$
on the time interval $\bar{I}_n = [t_{n-1},t_{n}]$.
Then, we obtain for each $i = 0, \dots, k-1$
\begin{equation} \label{lin_evol_In}
 \int_{I_n} \spV{d_t u_h(t), v} \psi_{n,i}(t) dt = \int_{I_n}
   \spV{F(t, u_h(t)), v} \psi_{n,i}(t) dt
 \qquad\forall\;  v\in W .
\end{equation}
By the definition of the weak time derivative we get for $u_h$ represented by
\eqref{poly_ansatz} the equation
\begin{equation*} \label{dtu_eq}
  \int_{I_n} \spV{d_t u_h(t), v} \psi_{n,i}(t) dt
  =
  \int_{I_n}\sum_{j=0}^{k} \spH{U_n^j, v} \phi_{n,j}'(t) \psi_{n,i}(t)\, dt
 \qquad\forall\;  v\in W .
\end{equation*}
We define the basis functions $\phi_{n,j}\in\P_{k}(\bar{I}_n)$ of \eqref{poly_ansatz}
via the reference transformation
\mbox{$\omega_n : \Ih \to \bar{I}_n$} where $\Ih := [-1,1]$
and
\begin{equation*} \label{ref_tra}
  t = \omega_n(\th) := \frac{t_{n-1}+t_n}{2} + \frac{h_n}{2} \th \in \bar{I}_n
  \qquad\forall\; \th\in\Ih , \; n=1,\ldots, N.
\end{equation*}
Let $\phih_j\in\P_{k}(\Ih)$, $j=0,\ldots, k$, be
suitable basis functions  satisfying the conditions
\begin{equation} \label{phi_lr}
  \phih_j(-1) = \delta_{0,j}, \qquad
  \phih_j(1) = \delta_{k,j},
\end{equation}
where $\delta_{k,j}$ denotes the usual Kronecker symbol.
Then, we define the basis functions on the original time interval $\bar{I}_n$
by
\begin{equation*} \label{phi_ref}
  \phi_{n,j}(t) := \phih_j(\th) \qquad\text{with}\qquad
  \th := \omega_n^{-1}(t) =
  \frac{2}{h_n} \left( t - \frac{t_n - t_{n-1}}{2}  \right) \in \Ih .
\end{equation*}
Similarly, we define the test basis functions $\psi_{n,i}$ by suitable reference
basis functions $\psih_i\in\P_{k-1}(\Ih)$, i.e.,
\begin{equation*} \label{psi_ref}
  \psi_{n,i}(t) := \psih_i( \omega_n^{-1}(t) )
  \qquad\forall\; t\in \bar{I}_n ,\; i = 0, \ldots , k-1 .
\end{equation*}
By the property \eqref{phi_lr}, the initial condition and the continuity
(with respect to time) of the discrete solution $u_h:I\to W$ is equivalent
to the conditions:
\begin{equation*} \label{U-left}
  U_1^{0} = u_0   \qquad\text{and}\qquad
  U_n^{0} = U_{n-1}^{k} \quad\forall\; n>2 .
\end{equation*}
We transform the integrals in \eqref{lin_evol_In} to the reference interval
$\Ih$ and obtain the following
system of equations for the ''coefficients'' $U_n^j\in W$, $j=1,\ldots,k$,
in the ansatz \eqref{poly_ansatz} :
\begin{equation} \label{lin_evol_sys}
  \sum_{j=0}^{k} \alpha_{i,j} \spH{U_n^j,v}
  =
  \frac{h_n}{2} \int_{\Ih} \spV{
  F\left( \omega_n(\th), \sum_{j=0}^{k} U^j_n \phih_j(\th) \right), v}
  \psih_i(\th) \,d\th
  \qquad\forall\; v\in W
\end{equation}
where $i=0,\ldots,k-1$,
$$
  \alpha_{i,j} := \int_{\Ih} {d_{\th} \phih_{j}}(\th) \psih_i(\th) \,d\th ,
$$
and the ''coefficient'' $U_n^0\in W$ is known.
We approximate the integral on the right hand side of \eqref{lin_evol_sys}
by the ($k+1$)-point {\em Gau{\ss}-Lobatto quadrature formula}:
\begin{equation*} \label{f_GL_approx}
    \int_{\Ih}  \spV{
  F\left( \omega_n(\th), \sum_{j=0}^{k} U^j_n \phih_j(\th) \right), v}
    \psih_i(\th) \,d\th
  \;\approx \;
  \sum_{\mu=0}^{k}
    \ww_{\mu}  \spV{
  F\left( \omega_n(\th_{\mu}), \sum_{j=0}^{k} U^j_n \phih_j(\th_{\mu}) \right), v}
     \psih_i(\th_{\mu}),
\end{equation*}
where $\ww_{\mu}$ are the weights and $\th_{\mu}\in[-1,1]$ are the integration points
with $\th_0=-1$ and $\th_k=1$.
Let us define the mapped Gau{\ss}-Lobatto points $t_{n,\mu}\in \bar{I}_n$ and the
coefficients $\beta_{i,\mu}$, $\gamma_{j,\mu}$  by
$$
   t_{n,\mu}:=\omega_n(\th_{\mu}) , \qquad
   \beta_{i,\mu} := \ww_{\mu}  \psih_i(\th_{\mu}) , \qquad
   \gamma_{j,\mu} := \phih_j(\th_{\mu}) .
$$
Then, the system \eqref{lin_evol_sys} is equivalent to the following system of
equations for the $k$ unknown ''coefficients'' $U_n^j\in W$, $j=1,\ldots,k$,
\begin{equation} \label{lin_sys}
  \sum_{j=0}^{k} \alpha_{i,j} \spH{U_n^j, v}
  =
  \frac{h_n}{2} \sum_{\mu=0}^{k} \beta_{i,\mu} \spV{
  F\left( t_{n,\mu}, \sum_{j=0}^{k} \gamma_{j,\mu} U^j_n \right), v}
  \qquad\forall\; v\in W.
\end{equation}
with the $k$ ''equations'' $i=0,\ldots,k-1$ where $U^0_n = U^k_{n-1}$ for
$n>1$ and $U^0_1= u_0$.

Once we have solved this system we enter the next time interval and set the
initial value of the new time interval $I_{n+1}$  to $U_{n+1}^0 := U_n^{k}$.
If the Gau{\ss}-Lobatto formula would be exact for the right hand side of
\eqref{lin_evol_sys}
this time marching process would solve the global time discretization
\eqref{cGP_form} exactly.
Since in general there is an integration error
we call the time marching process corresponding to  \eqref{lin_sys} simply the
''cGP(k)-{\em method}''.

In principle, we have to solve a coupled system for the $U_n^j\in W$ which
could be very expensive.
However, by a clever choice of the functions $\phih_j$ and $\psih_i$
it is possible to uncouple the system to a large extend.
In the following, we will discuss this issue for the special methods
cGP(1), cGP(2) and for the general method cGP($k$), $k\ge 3$.
In all cases,
we choose the basis functions $\phih_j\in\P_k(\Ih)$
as the Lagrange basis functions with respect to the Gau{\ss}-Lobatto
points $\th_{\mu}$, i.e.,
\begin{equation*} \label{phih_choice}
  \phih_j(\th_{\mu}) =  \delta_{j,\mu}
  \qquad\forall\; j, \mu\in \{0, \ldots, k \}.
\end{equation*}
Then, the method \eqref{lin_sys} reduces to
\begin{equation*} \label{lin_sys_1}
  \sum_{j=0}^{k} \alpha_{i,j} \spH{U_n^j, v}
  =
  \frac{h_n}{2} \sum_{j=0}^{k} \beta_{i,j} \spV{
  F\left( t_{n,j},  U^{j}_n \right), v}
  \qquad\forall\; v\in W,\;  i = 0, \ldots, k-1,
\end{equation*}
and by the choice of the test basis functions $\psih_i\in\P_{k-1}(\Ih)$
we try to
get suitable values for the coefficients $\alpha_{i,j}$ and
$\beta_{i,j}$.
In the following, we will use the following abbreviation and
assumption:
\begin{equation} \label{F_j}
     F_n^j(U_n^j) := F(t_{n,j}, U_n^j) \in H'
     \qquad\forall\; j = 0, \ldots, k, \; n = 1, \ldots, N.
\end{equation}
\subsubsection{The cGP(1) method%
\label{cGP1_method}
}
We use the 2-point Gau{\ss}-Lobatto formula (trapezoidal rule) with
$\ww_0 = \ww_1 = 1$ and $\th_0=-1$, $\th_1=1$.
The only test function
$\psih_0$ is chosen as $\psih_0(\th)=1$.
Then, we obtain
$$
 \alpha_{0,0}=-1, \quad  \alpha_{0,1}=1, \quad
 \beta_{0,0} = \beta_{0,1} = 1 .
$$
Using the notation $U^{n-1} := u_h(t_{n-1}) = U_n^0$ and
$U^{n} := u_h(t_{n}) = U_n^1$, we obtain the following equation for
the ''unknown'' $U^n\in W$ :
\begin{equation*} \label{cGP1_eq}
 \spH{U^n,v} - \spH{U^{n-1},v}
 =
 \frac{h_n}{2} \left\{
   \spV{ F(t_{n-1}, U^{n-1}) + F(t_{n}, U^n), v }
 \right\}
\end{equation*}
for all $v\in W$
which is the well-known {\em Crank-Nicolson method}.
In operator notation it can be written in the equivalent form:
\begin{equation*} \label{cGP1_M_form}
  U^n
 =  U^{n-1} +
 \frac{h_n}{2} \Id^{-1}\left\{
   F(t_{n-1}, U^{n-1}) + F(t_{n}, U^n)
 \right\}.
\end{equation*}

\subsubsection{The cGP(2) method%
\label{cGP2_method}
}
We use the 3-point Gau{\ss}-Lobatto formula (Simpson rule) with
$\ww_0 = \ww_2 = 1/3$, $\ww_1 = 4/3$ and $\th_0=-1$, $\th_1=0$, $\th_2=1$.
For  the test functions  $\psih_i\in\P_1(\Ih)$, we choose
$$
 \psih_0(\th) = - \frac{3}{4} \th, \quad
 \psih_1(\th) = 1 .
$$
Then, we get
$$
 (\alpha_{i,j}) =
     \begin{pmatrix}
       -1/2  &  1  & -1/2  \\
      -1  &     0  &  1
     \end{pmatrix},
 \quad
 (\beta_{i,j}) =
     \begin{pmatrix}
       1/4    &    0  & -1/4  \\
       1/3    &  4/3  &  1/3
     \end{pmatrix}
 \quad
$$
and the assumption~\eqref{F_j},
the system to compute the ''unknowns'' $U_n^1, U_n^2 \in W$ from
the known $U_n^0 = U_{n-1}^2$ reads:
  \begin{eqnarray}
     \label{cGP2_1}
     U_n^1 &=&
         \frac{1}{2} U_n^0 + \frac{1}{2} U_n^2
        +
        \frac{h_n}{8} M^{-1}
          \left\{ F_n^0(U_n^0) - F_n^2(U_n^2) \right\}
    \\[1mm]
     \label{cGP2_2}
    U_n^2
     &=& U_n^0 + \frac{h_n}{6}  M^{-1}
        \left\{ F_n^0(U_n^0) + 4 F_n^1(U_n^1) + F_n^2(U_n^2) \right\} .
  \end{eqnarray}
Let us denote the value for $U_n^1$ computed from \eqref{cGP2_1} and depending
on $U_n^2$ by $U_n^1 = G_n^1( U_n^2 )$
where $G_n^1: W \to W$ in general is a nonlinear operator.
We substitute this in the equation  \eqref{cGP2_2} and get, for the
unknown $U_n^2\in W$, the following fixed point equation :
\begin{equation*} \label{cGP2_fixpoint}
  U_n^2 = G_n^2 (U_n^2)
  := U_n^0 + \frac{h_n}{6} \Id^{-1}\left\{
    F_n^0(U_n^0)  + 4 F_n^1( G_n^1(U_n^2) ) + F_n^2(U_n^2)
  \right\}
\end{equation*}
The mapping $G_n^2 : W \to W$ is a contraction if the time step size $\tn$
is sufficiently small.

\subsection{Gauss implicit Runge-Kutta methods}
For the general autonomous first order differential equations 
\begin{equation} \label{aut_de}
y'(t)=f(y(t)),
\end{equation}
where for system \eqref{hamiltonian}, we choose $
y = \left( \begin{array}{cc}
p  \\
q \end{array} \right)$  and $f(y) = \left( \begin{array}{cc}
-\triangledown_q V(q)  \\
-\triangledown_p T(q) \end{array} \right).$ Runge-Kutta methods are defined as

\begin{equation*}
y_{n+1}=y_{n}+h\sum _{i=1}^{s} b_{i}f(Y_{i})
\end{equation*}
and
\begin{equation*}
Y_{i}=y_{n}+h\sum _{i=1}^{s} a_{ij}f(Y_{j})
\end{equation*}

where the coefficients $a_{ij}$, $b_i$ and stage $s$ determine the method. The Gauss methods have the highest possible order $r=2s$ and are symplectic and symmetric.
We exclusively consider $s=2$, fourth order method for a fair comparison.

\subsection{General linear methods}
General linear methods provide numerical solutions of initial value problems of the form \eqref{aut_de} 
A general linear method is of the form,
\begin{align*}
Y &= h(A\otimes I)f(Y)+ (U\otimes I)y^{[n-1]},\\
y^{[n]}&= h(B\otimes I)f(Y)+ (V\otimes I)y^{[n-1]}.
\end{align*}
where $A\otimes I$ is the Kronecker product of the matrix $A$ and the identity matrix $I$ and $h$ represents the step size. The $s-$component vector $Y$ are the stages and $f(Y)$ are the stage derivatives. The vector $y^{[n-1]}$ with $r-$components is an input at the beginning of a step and results in output approximation $y^{[n]}$.
With a slight abuse of notation, we can write,
\begin{align*}\label{general_lin}
Y &= hAf(Y)+ Uy^{[n-1]},\nonumber \\
y^{[n]} &= hBf(Y)+ Vy^{[n-1]}.
\end{align*}
The matrices $A$, $U$, $V$ and $B$ represent a particular general linear method and are generally displayed as,
\begin{equation*}\label{matrix_glm}
 \left[ {\begin{array}{c|c}
 A & U \\ \hline
B & V 
 \end{array} } \right].
\end{equation*}

A fourth order symmetric $G$-symplectic general linear method is constructed with four stages $(s=4)$ and three input values $(r=3)$. The coefficeints of 
the method are given in \cite{butcher_yousaf}.

\section{Numerical Experiments}
We performed numerical comparisons of the continuous Galerkin Petrov time discretization scheme , general linear method and implicit Gauss R-K method all having the same order four, for some Hamiltonian systems including
simple harmonic oscillator, Kepler's problem with different eccentricities and molecular dynamical problems. Throughout the comparison, continuous Galerkin
Petrov scheme is denoted by acronym cGP(2), while general linear method and implicit Gauss R-K method are represented by the acronym glm4 and irk4 respectively.
The emphasis in our comparison is on the accuracy of solution, including the phase information, enrgy conservation and CPU time using above discussed methods.
 For each method and problem, we used different stepsizes and several intervals of integration. Stepsizes
were chosen as a compromise between having small truncation error and performing efficient integration on each step.
The accuracy of the solution was measured by the $L_2$ norm of the absolute global error in the position and velocity coordinates and is denoted by $E_g(t)$.
The relative error in Hamiltonian is defined as $$E_e(t)=\frac{E(t)-E(0)}{E(0)}.$$ Growth of global error is measured for first two problems
as their exact solution exists, while relative error in Hamiltonian $E_e(t)$ is calcultaed for all problems. We also measured computational effort using the CPU time.
All the comparisons are done on the same machine and are optimized using MATLAB.

\subsection*{Simple Harmonic Oscillator}
As an example of simple harmonic oscillator a mass spring system having kinetic energy ${p^{2}}/({2m})$, where $p=mv$ is the momentum
of the system and potential energy $\frac{1}{2}kq^{2}$. Where $q$ is distance from the equilibrium, $m$ is the mass of the body which is
attached to spring and $k$ is constant of proportionality often called as spring constant. Here the Hamiltonian
is the total energy of the system and has one degree of freedom
\begin{equation*}
 H(q,p)=\frac{1}{2}kq^{2} + \frac{p^{2}}{2m}.
\end{equation*}
 The equations of motion from the Hamiltonian are
\begin{equation*}
 q'=\frac{\partial H(q,p)}{\partial p}=p,\;\;\;\;\;\;\;\;\:\:\:\:p'=-\frac{\partial H(q,p)}{\partial q}=-q.
\end{equation*}

\begin{figure}[H]
\begin{center}
\hspace{-2.5in}
\includegraphics[width=5in,height=2in]{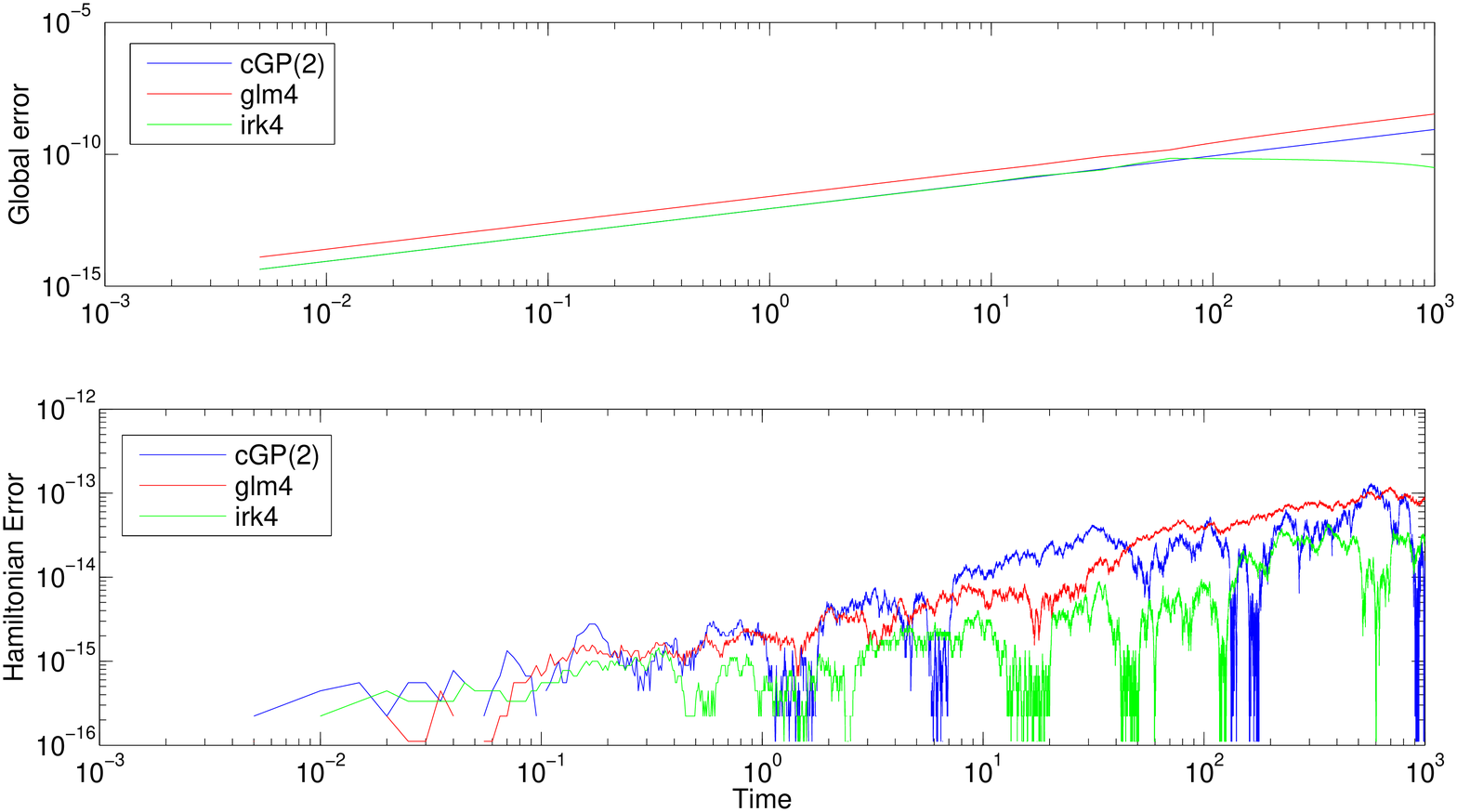}
\vspace{0.1in}
\caption{The growth of global error and relative error in Hamiltonian for Simple harmonic oscillator using stepsize $h=0.005$.}\label{fig:sho_ge}
\end{center}
\end{figure}

We compared the problem using different stepsizes of $h=0.005, 0.01, 0.025,$ and $0.05$. 
Figure \ref{fig:sho_ge} gives the log-log graph for time versus global error $E_g(t)$ and relative error in Hamiltonian $E_e(t)$ using stepsize $h=0.005$ 
for the time interval [0, 1000]. We found almost the same behavior of error growth for position and Hamiltonian using the rest of stepsizes.
In Figure \ref{fig:sho_ge}, the top plot gives the growth of global error and is approximately same for all tested methods, irk4 and cGP(2) having the least error
while glm4 with slightly bigger error. In bottom plot of figure \ref{fig:sho_ge}, the error in Hamiltonian is conserved by the methods. 
We also calculated the error growth according to Brouwer's law~\cite{bl}, our calculation shows that the exponent of time is 1 and 0.6 for $E_g(t)$ and $E_e(t)$
respectively, closed to its expected value.
Table \ref{sho_table} gives the cost of integration for simple harmonic oscillator using all stepsizes. The table lists the stepsizes, maximum of global error,
maximum of Hamiltonian error and CPU time. We observe from the Table \ref{sho_table} that cGP(2) used the least CPU time and also having the least value for
maximum of global error except for $h=0.005$, where irk4 having the least end point global error, may be because of entering in a dip also depicted in Figure \ref{fig:sho_ge}.
The methods irk4 and glm4 are using eight and sixteen times more CPU time than cGP(2) giving similar accuracy for $h=0.005$.

 \begin{table}
 \begin{center}
{\footnotesize
 \begin{tabular} {l|cccc}

\multicolumn{1}{l|}{Method} & \multicolumn{1}{c}{stepsize $(h)$} & 
\multicolumn{1}{c}{Max. of Global} &
\multicolumn{1}{c}{Max. of Hamiltonian} &
\multicolumn{1}{c}{CPU Time (sec.)} \\
\multicolumn{1}{l|}{} & {}& \multicolumn{1}{c}{Error} &
\multicolumn{1}{c}{Error} &
\multicolumn{1}{c}{} \\
\hline
\\
cGP(2)  & $0.05$ & $8.67\times10^{-6}$ & $4.08\times10^{-14}$& $2.8$\\
cGP(2)  & $0.025$ & $5.42\times10^{-7}$ & $1.07\times10^{-13}$& $5.3$\\
cGP(2)  & $0.01$ & $1.38\times10^{-8}$ & $1.12\times10^{-13}$& $12.1$\\
cGP(2)  & $0.005$ & $8.68\times10^{-10}$ & $1.31\times10^{-13}$& $23.9$\\
 \\
glm4  & $0.05$ & $2.51\times10^{-5}$ & $4.67\times10^{-11}$& $3.2$\\
glm4  & $0.025$ & $1.57\times10^{-6}$ & $7.57\times10^{-13}$& $11.6$\\
glm4  & $0.01$ & $4.09\times10^{-8}$ & $3.1\times10^{-15}$& $78.6$\\
glm4  & $0.005$ & $3.34\times10^{-9}$ & $1.14\times10^{-15}$& $395.3$\\
\\
irk4  & $0.05$ & $8.67\times10^{-6}$ & $6.43\times10^{-15}$& $5.7$\\
irk4  & $0.025$ & $5.41\times10^{-7}$ & $2.51\times10^{-14}$& $11.5$\\
irk4  & $0.01$ & $1.31\times10^{-8}$ & $1.97\times10^{-14}$& $45.0$\\
irk4  & $0.005$ & $6.98\times10^{-11}$ & $3.68\times10^{-14}$& $191.8$\\

 \end{tabular}
}\footnotesize
\caption{Maximum of global error, Hamiltonian error and CPU time for simple harmonic oscillator.}
 \label{sho_table}
 \end{center}
\end{table}

\subsubsection*{Kepler's Problem}
Kepler's problem is two body orbital problem in which the bodies are moving under their mutual gravitational forces. We can assume that one
body is fixed at the origin and the second body is located in the plane with coordinates $(q_{1},q_{2}).$  The solution of this problem is used in many
important applications which includes the determination of orbits for new asteroids and the measurement of orbits for the two primary bodies in a restricted three body problem. The Hamiltonian of the system
can be written in separable form as \cite{hairer_gni}
\begin{equation*}
 H(q,p)=\frac{1}{2} (p_{1}^{2}+p_{2}^{2}) - \frac{1}{\sqrt{q_{1}^{2}+q_{2}^{2}}}
\end{equation*}
This can be written as $H = T + V$, where $T = (p_{1}^{2}+p_{2}^{2})/{2}$ and $V = -1/{\sqrt{q_{1}^{2}+q_{2}^{2}}}$ are kinetic and
potential energy of the system respectively. As like the previous problem, this system is also autonomous so the Hamiltonian $H$ is
a conserved quantity. \newline
The equations of motion are
\begin{equation*}
 q'_{1}=p_{1},\:\:\:\:\:\:\;\;\;\:\;\:q'_{2}=p_{2}
\end{equation*}
\begin{equation}\label{kepler}
 p'_{1}=q''_{1}=-\frac{q_{1}}{(q_{1}^2+q_{2}^2)^\frac{3}{2}}
\end{equation}
\begin{equation*}
 p'_{2}=q''_{2}=-\frac{q_{2}}{(q_{1}^2+q_{2}^2)^\frac{3}{2}}
\end{equation*}
with the initial conditions
\begin{equation*}
 q_{1}(0)=1-e,\;\;\;q_{2}(0)=0,\;\;\;q'_{1}(0)=0,\;\;\;q_{2}(0)=\sqrt{\frac{1+e}{1-e}}
\end{equation*}
where $e$ is eccentricity $0 \leq e < 1$. The exact solution of the above equations (\ref{kepler}) is
\begin{equation*}
y_{1} = \cos(E) - e,\;\;\;\;\;y_{2} = \sqrt{1-e^{2}}\sin(E),
\end{equation*}
and
\begin{equation*}
y_{1}' = -\sin(E)(1-e\cos(E))^{-1},\;\;\;\;\;y_{2}' = \sqrt{(1-e^{2})}\cos(E)(1-e\cos(E))^{-1},
\end{equation*}
where the eccentric anomaly $E$ satisfies Kepler's equation $t=E-e\sin(E)$. Since Kepler's equation is implicit in $E$, the equation is usually
solved using a non-linear equation solver, although useful analytical approximations can be found for smaller eccentricity.

\begin{figure}[H]
\begin{center}
\hspace{-2.5in}
\includegraphics[width=5in,height=2in]{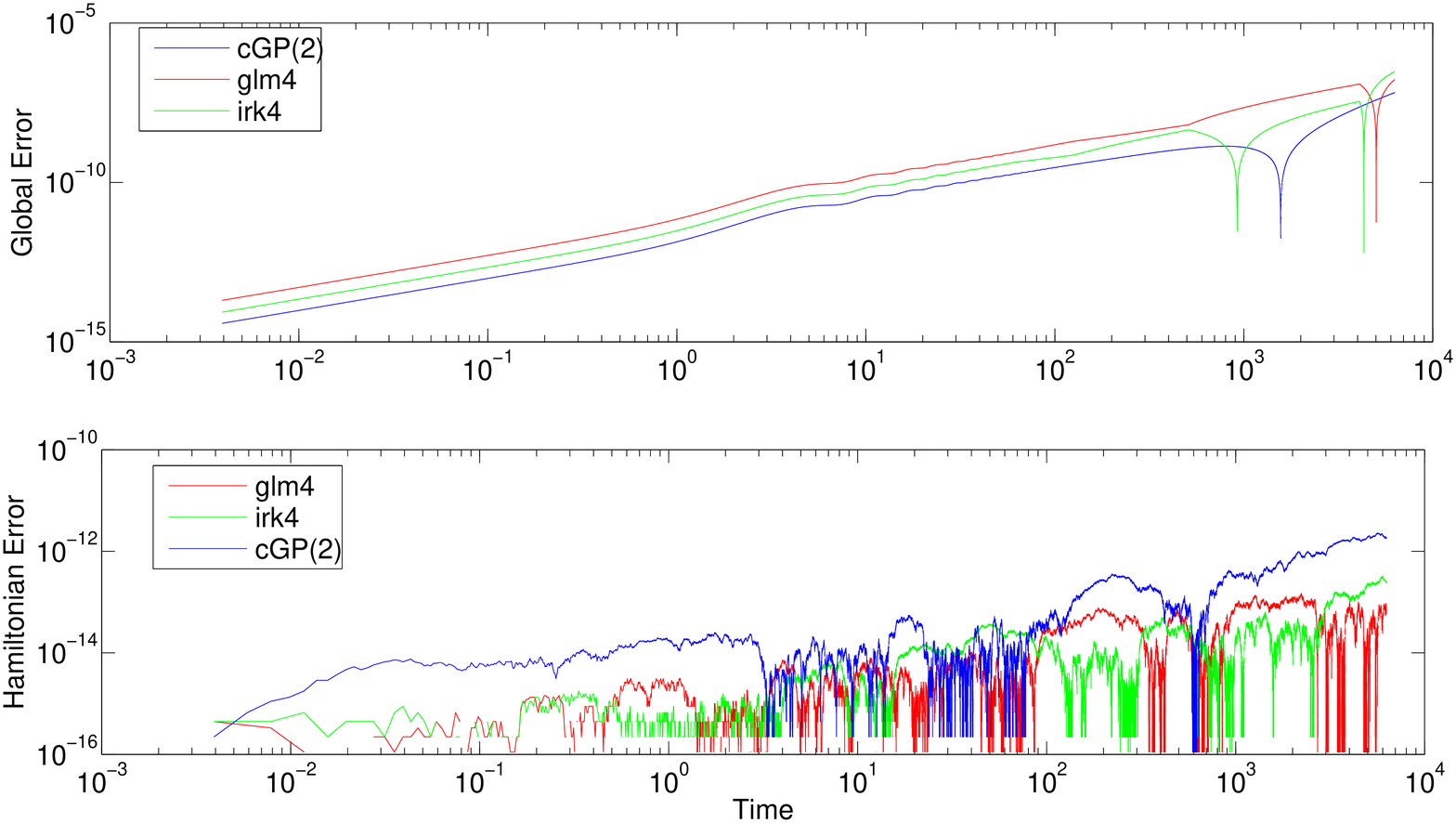}
\caption{The growth of global error and relative error in Hamiltonian for kepler's problem with $e=0$ using stepsize $2\pi/6400$ for $10^3$ periods.}\label{fig:kep0_ge}
\end{center}
\end{figure}

\begin{figure}[H]
\begin{center}
\hspace{-2.5in}
\includegraphics[width=5in,height=2in]{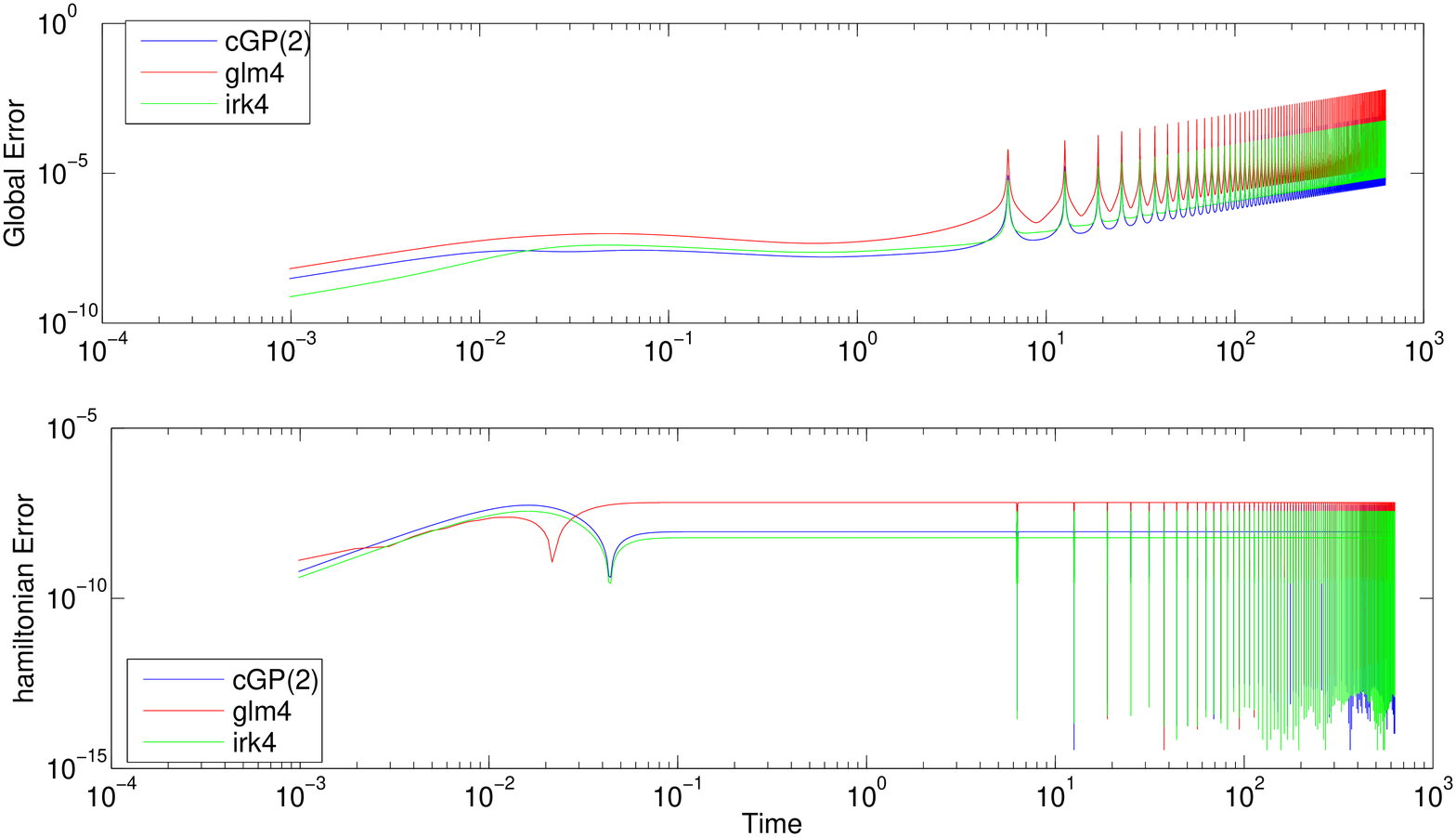}
\caption{The growth of global error and relative error in Hamiltonian for kepler's problem with $e=0.9$ using stepsize $2\pi/1600$ for $10^2$ periods.}\label{fig:kep0p9_ge}
\end{center}
\end{figure}

The integrations are performed for Kepler's problem with different eccentricities $e=0, 0.5~ \mathtt{and}~ 0.9$. The integration is done for
$1000$ periods for $e=0$ and 100 periods for $e=0.5 ~ \mathtt{and}~ 0.9$. For each method, we measured $E_g(t)$ and $E_e(t)$ throughout the interval of integration. A variety of different stepsizes
are used to analyze the behaviour of error growth. We used the stepsizes of $h=\frac{2\pi}{400}$, $h=\frac{2\pi}{800}$ $h=\frac{2\pi}{1600}$, 
$h=\frac{2\pi}{3200}$ and $h=\frac{2\pi}{6400}$ for eccentricities $e=0,0.5 ~\mathtt{and}~ 0.9$. 
A log-log plot of time against error is given for Kepler's problem in Figures \ref{fig:kep0_ge} and \ref{fig:kep0p9_ge} using eccentricities 0 and 0.9 respectively.
Growth of errors in both quantities behave in the same manner as for e=0.5. It is seen that the global error growth is approximately linear for cGP(2), irk4 and glm4, i. e., growing as $t^{0.9}$ (see figures \ref{fig:kep0_ge} and \ref{fig:kep0p9_ge}).
The error in Hamiltonian remains conserved  for cGP(2), irk4 and glm4 for the intervals of integration. Our calculation shows that for $E_e(t)$ grows as $t^{0.6}$, showing a good agreement to its expected value.
The cGP(2) exhibits a smaller error even the problem becomes more eccentricitic (see Figures \ref{fig:kep0_ge} and \ref{fig:kep0p9_ge}).

We also measured the cost of integration for Kepler's problem using all stepsizes for all three eccentricities. Tables \ref{kepe0_table}, \ref{kep0p5_table} 
and  \ref{kep0p9_table} lists the stepsizes, maximum of global error, maximum of Hamiltonian error and CPU time for $e=0, e=0.5 ~ \mathtt{and}~ 0.9$ respectively.
We observe from the information depicted in tables, that cGP(2) used the least CPU time and also having the least value for maximum of global error for all the stepsizes.
For e=0, using the least stepsize i.e $h=\frac{2\pi}{6400}$, irk4 and glm4 used 506 and 26 times more CPU time than cGP(2). While for e=0.5 and 0.9,  irk4 and glm4 used nearly 55 and 24 times more CPU time than cGP(2).

 \begin{table}
 \begin{center}
{\footnotesize
 \begin{tabular} {l|cccc}

\multicolumn{1}{l|}{Method} & \multicolumn{1}{c}{stepsize $(h)$} & 
\multicolumn{1}{c}{Max. of Global} &
\multicolumn{1}{c}{Max. of Hamiltonian} &
\multicolumn{1}{c}{CPU Time (sec.)} \\
\multicolumn{1}{l|}{} & {}& \multicolumn{1}{c}{Error} &
\multicolumn{1}{c}{Error} &
\multicolumn{1}{c}{} \\
\hline
\\
cGP(2)  & $2\pi/400$ & $4.88\times10^{-6}$ & $4.07\times10^{-13}$& $93.4$\\
cGP(2)  & $2\pi/800$ & $2.54\times10^{-7}$ & $7.54\times10^{-12}$& $192.7$\\
cGP(2)  & $2\pi/1600$ & $6.56\times10^{-8}$ & $1.28\times10^{-11}$& $378$\\
cGP(2)  & $2\pi/3200$ & $1.38\times10^{-8}$ & $2.28\times10^{-12}$& $760.5$\\
cGP(2)  & $2\pi/6400$ & $4.54\times10^{-9}$ & $1.49\times10^{-12}$& $1487$\\
 \\
glm4  & $2\pi/400$ & $2.36\times10^{-5}$ & $2.35\times10^{-14}$& $3464$\\
glm4  & $2\pi/800$ & $1.56\times10^{-6}$ & $3.06\times10^{-14}$& $12172$\\
glm4  & $2\pi/1600$ & $1.66\times10^{-7}$ & $9.39\times10^{-14}$& $47724$\\
glm4  & $2\pi/3200$ & $8.59\times10^{-8}$ & $5.32\times10^{-13}$& $180716$\\
glm4  & $2\pi/6400$ & $1.05\times10^{-8}$ & $9.98\times10^{-12}$& $752864$\\
\\
irk4  & $2\pi/400$ & $1.04\times10^{-5}$ & $4.72\times10^{-14}$& $1598$\\
irk4  & $2\pi/800$ & $5.68\times10^{-7}$ & $2.79\times10^{-14}$& $6348$\\
irk4  & $2\pi/1600$ & $2.98\times10^{-7}$ & $3.28\times10^{-14}$& $23609$\\
irk4  & $2\pi/3200$ & $6.58\times10^{-8}$ & $5.17\times10^{-13}$& $93215$\\
irk4  & $2\pi/6400$ & $1.27\times10^{-8}$ & $7.18\times10^{-12}$& $39460$\\

 \end{tabular}
}\footnotesize
\caption{Maximum of global error, Hamiltonian error and CPU time for Kepler's Problem with $e=0$ for $10^3$ periods.}
 \label{kepe0_table}
 \end{center}
\end{table}

 \begin{table}
 \begin{center}
{\footnotesize
 \begin{tabular} {l|cccc}

\multicolumn{1}{l|}{Method} & \multicolumn{1}{c}{stepsize $(h)$} & 
\multicolumn{1}{c}{Max. of Global} &
\multicolumn{1}{c}{Max. of Hamiltonian} &
\multicolumn{1}{c}{CPU Time (sec.)} \\
\multicolumn{1}{l|}{} & {}& \multicolumn{1}{c}{Error} &
\multicolumn{1}{c}{Error} &
\multicolumn{1}{c}{} \\
\hline
\\
cGP(2)  & $2\pi/400$ & $1.06\times10^{-4}$ & $2.83\times10^{-8}$& $10.2$\\
cGP(2)  & $2\pi/800$ & $6.63\times10^{-6}$ & $1.77\times10^{-9}$& $19.9$\\
cGP(2)  & $2\pi/1600$ & $4.15\times10^{-7}$ & $1.11\times10^{-10}$& $39.8$\\
cGP(2)  & $2\pi/3200$ & $2.54\times10^{-8}$ & $7.01\times10^{-12}$& $79.2$\\
cGP(2)  & $2\pi/6400$ & $2.9\times10^{-9}$ & $1.02\times10^{-12}$& $157.5$\\
 \\
glm4  & $2\pi/400$ & $2.1\times10^{-4}$ & $5.97\times10^{-8}$& $43.8$\\
glm4  & $2\pi/800$ & $1.31\times10^{-5}$ & $3.74\times10^{-9}$& $182.3$\\
glm4  & $2\pi/1600$ & $8.2\times10^{-7}$ & $2.34\times10^{-10}$& $688.2$\\
glm4  & $2\pi/3200$ & $3.31\times10^{-8}$ & $1.46\times10^{-11}$& $2366$\\
glm4  & $2\pi/6400$ & $2.16\times10^{-8}$ & $1.38\times10^{-12}$& $8789$\\
\\
irk4  & $2\pi/400$ & $8.84\times10^{-5}$ & $1.89\times10^{-8}$& $22.1$\\
irk4  & $2\pi/800$ & $5.53\times10^{-6}$ & $1.18\times10^{-9}$& $68.5$\\
irk4  & $2\pi/1600$ & $3.43\times10^{-7}$ & $4.71\times10^{-11}$& $244$\\
irk4  & $2\pi/3200$ & $1.04\times10^{-8}$ & $4.66\times10^{-12}$& $955$\\
irk4  & $2\pi/6400$ & $2.51\times10^{-8}$ & $3.13\times10^{-13}$& $3830$\\

 \end{tabular}
}\footnotesize
\caption{Maximum of global error, Hamiltonian error and CPU time for Kepler's Problem with $e=0.5$ for $10^2$ periods.}
 \label{kep0p5_table}
 \end{center}
\end{table}

 \begin{table}
 \begin{center}
{\footnotesize
 \begin{tabular} {l|cccc}

\multicolumn{1}{l|}{Method} & \multicolumn{1}{c}{stepsize $(h)$} & 
\multicolumn{1}{c}{Max. of Global} &
\multicolumn{1}{c}{Max. of Hamiltonian} &
\multicolumn{1}{c}{CPU Time (sec.)} \\
\multicolumn{1}{l|}{} & {}& \multicolumn{1}{c}{Error} &
\multicolumn{1}{c}{Error} &
\multicolumn{1}{c}{} \\
\hline
\\
cGP(2)  & $2\pi/400$ & $4.87$ & $5.23\times10^{-3}$& $10.4$\\
cGP(2)  & $2\pi/800$ & $1.68$ & $2.43\times10^{-4}$& $20.2$\\
cGP(2)  & $2\pi/1600$ & $1.92\times10^{-1}$ & $1.42\times10^{-5}$& $41.1$\\
cGP(2)  & $2\pi/3200$ & $1.3\times10^{-2}$ & $8.74\times10^{-7}$& $80.6$\\
cGP(2)  & $2\pi/6400$ & $8.41\times10^{-6}$ & $5.43\times10^{-8}$& $162.1$\\
 \\
glm4  & $2\pi/400$ & $111.8$ & $6.5\times10^{-4}$& $44.1$\\
glm4  & $2\pi/800$ & $4.54$ & $2.23\times10^{-4}$& $172$\\
glm4  & $2\pi/1600$ & $1.46$ & $1.62\times10^{-5}$& $676$\\
glm4  & $2\pi/3200$ & $9.77\times10^{-2}$ & $1.04\times10^{-6}$& $2482$\\
glm4  & $2\pi/6400$ & $6.14\times10^{-3}$ & $6.53\times10^{-8}$& $8402$\\
\\
irk4  & $2\pi/400$ & $4.54$ & $1.73\times10^{-5}$& $22.4$\\
irk4  & $2\pi/800$ & $1.27$ & $1.36\times10^{-5}$& $67.8$\\
irk4  & $2\pi/1600$ & $1.32\times10^{-1}$ & $1.38\times10^{-6}$& $257.5$\\
irk4  & $2\pi/3200$ & $9.01\times10^{-3}$ & $9.45\times10^{-8}$& $1086$\\
irk4  & $2\pi/6400$ & $5.73\times10^{-4}$ & $6.03\times10^{-9}$& $3786$\\

 \end{tabular}
}\footnotesize
\caption{Maximum of global error, Hamiltonian error and CPU time for Kepler's Problem with $e=0.9$ for $10^2$ periods.}
 \label{kep0p9_table}
 \end{center}
\end{table}

\subsection*{Molecular Dynamical Problem}
We consider the interaction of seven Argon atoms in two dimension, where one of the atom is centered by six atoms which are symmetrically
arranged \cite{bies_skeel}. The Hamiltonian for the molecular dynamics is written as \cite{hairer_gni}
\begin{equation*}
H(q,p)=\frac{1}{2} \sum_{i=1}^{7} \frac{1}{m_{i}}p_{i}^{T}p_{i} + \sum_{i=2}^{7}\sum_{j=1}^{i-1} V_{ij}\Vert q_{i}-q_{j}\Vert
\end{equation*}
where $V_{ij}(r)$ are potential functions. Here $q_{i}$ and $p_{i}$ are positions and generelized momenta for the atoms. And $m_{i}$ denotes
the atomic mass of the $i$th atom.
\begin{equation*}
V_{ij}(r)=4\varepsilon_{ij}\left( \left( \frac{\sigma_{ij}}{r}\right)^{12} - \left( \frac{\sigma_{ij}}{r}\right)^6 \right).
\end{equation*}
The equations of motion for the frozen Argon crystals are given as
\begin{equation*}
q''_{i}(t)=\frac{24\varepsilon\sigma^6}{m_{i}} \sum_{j=1,j\neq i}^{7}\left[ \frac{(q_{j}-q_{i})}{\Vert q_{j}-q_{i}\Vert _{2}^{8}}-2\sigma^6\frac{(q_{j}-q_{i})}{\Vert q_{j}-q_{i}\Vert _{2}^{14}}\right ],\hspace{0.2in} i=1,...,7,
\end{equation*}
where $r=\sigma_{ij} \sqrt[6]{2}$, $m_{i}=66.34\times 10^{-27}\mathtt{[kg]}$, $\sigma_{ij}=\sigma=0.341\mathtt{[nm]}$ and $\varepsilon=1.654028284\times10^{-21}\mathtt{[J]}$.
Initial positions and initial velocities are taken in [nm] and [nm/sec] respectively \cite{hairer_gni}.

\begin{figure}[H]
\begin{center}
\hspace{-2.5in}
\includegraphics[width=5in,height=2in]{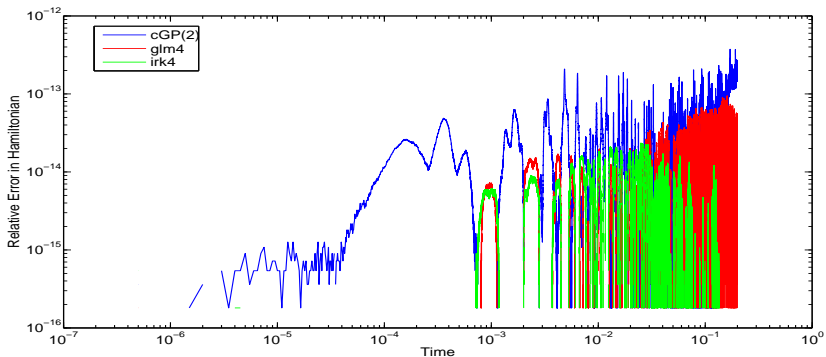}
\caption{The growth of relative error in Hamiltonian using  $h=0.5\times 10^{-6}\mathtt{[fsec]}$ for molecular dynamical problem over an interval of $2 \times 10^5$ [fsec].}\label{fig:mole_reha}
\end{center}
\end{figure}

In molecular dynamics, since much ineterst is emphasized on macroscopic quantities like Hamiltonian. So we also discussed only the energy
conservation of atoms over an interval of length $2 \times 10^5$ [fsec] ($1\mathtt{fsec}=10^{-6}$). The experiments are done using the
stepsizes of 0.5 fsec, 1 fsec, 2 fsec and 4 fsec. The graphical results are only shown for $h=0.5\times 10^{-6}\mathtt{[fsec]}$ as the error growth using other stepsizes
was approximately same. Figure \ref{fig:mole_reha} shows that the tested methods conserve the value of Hamiltonian $H$ even though the conservation is of 
highly oscillatory, while the error in Hamiltonian for cGP(2) grows as $t^{0.7}$. On the other hand, for irk4 and glm4 the exponent of time is 0.59 and 0.61 respectively.
Table \ref{mole_table} gives the cost of integration for molecular dynamical problem using all stepsizes. The table lists the stepsizes, maximum of Hamiltonian
error and CPU time. It is observed from the Table \ref{mole_table} that cGP(2) used the least CPU time for all the stepsizes used but exhibiting slightly 
big maximum of Hamiltonian error. The methods irk4 and glm4 having almost the same error growth for the integrated interval.

 \begin{table}
 \begin{center}
{\footnotesize
 \begin{tabular} {l|ccc}

\multicolumn{1}{l|}{Method} & \multicolumn{1}{c}{stepsize $(h)$} & 
\multicolumn{1}{c}{Max. of Global} &
\multicolumn{1}{c}{CPU Time (sec.)} \\
\multicolumn{1}{l|}{} & \multicolumn{1}{c}{} &
\multicolumn{1}{c}{Error} &
\multicolumn{1}{c}{} \\
\hline
\\
cGP(2)  & $4$ & $1.24\times10^{-9}$ &  $658$\\
cGP(2)  & $2$ & $6.15\times10^{-11}$ &  $1323$\\
cGP(2)  & $1$ & $7.05\times10^{-12}$ &  $2680$\\
cGP(2)  & $0.5$ & $3.73\times10^{-13}$ &  $5821$\\
 \\
glm4  & $4$ & $9.24\times10^{-11}$ &  $1772$\\
glm4  & $2$ & $5.58\times10^{-12}$ &  $4300$\\
glm4  & $1$ & $3.53\times10^{-13}$ &  $11367$\\
glm4  & $0.5$ & $9.41\times10^{-14}$ &  $34591$\\
\\
irk4  & $4$ & $3.72\times10^{-11}$ &  $2060$\\
irk4  & $2$ & $2.32\times10^{-12}$ &  $4272$\\
irk4  & $1$ & $1.59\times10^{-13}$ &  $10020$\\
irk4  & $0.5$ & $2.34\times10^{-14}$ &  $23596$\\

 \end{tabular}
}\footnotesize
\caption{Maximum of Hamiltonian error and CPU time for molecular dynamical problem for $2 \times 10^5$ [fsec].}
 \label{mole_table}
 \end{center}
\end{table}

\section{Summary}
We implemented and analyzed the cGP(2) for Hamiltonian systems such as harmonic oscillator, Kepler's problem and molecular dynamical problem. 
The obtained results are also compared with symplectic methods irk4 and glm4. It is shown that the cGP(2) method conserves the hamiltonian as other tested
symplectic methods do. Moreover, giving the efficiency approximately same as other methods yield, cGP(2) uses marginally less CPU time than compared methods.


\end{document}